\documentclass[reqno]{amsart}
\usepackage{amsmath,amssymb}
\usepackage{booktabs,mathrsfs,enumerate}
\usepackage{hyperref,framed}
\usepackage[hyperpageref]{backref}
\hypersetup{pdftitle={Semisimple representations of alternating cyclotomic Hecke algebras},
  pdfauthor={Clinton Boys and Andrew Mathas},
  colorlinks = true,
  linkcolor =blue,
  anchorcolor = red,
  citecolor = blue,
  urlcolor = blue
}

\renewcommand*{\backref}[1]{}
\renewcommand*{\backrefalt}[4]{%
  \ifcase #1 No citations.
  \or [Page #2.]%
  \else [Pages #2.]%
  \fi%
}

\usepackage{todonotes}

\usepackage{aliascnt}
\def\NewTheorem#1{%
  \newaliascnt{#1}{equation}
  \newtheorem{#1}[#1]{#1}
  \aliascntresetthe{#1}
  \expandafter\def\csname #1autorefname\endcsname{#1}
}
\def\equationautorefname~#1\null{(#1)\null}

\theoremstyle{plain}
\newcounter{mainTheorem}

\newcounter{mainCorollary}
\numberwithin{mainCorollary}{mainTheorem}

\newtheorem{Corollary*}[mainCorollary]{Corollary}

\swapnumbers
\numberwithin{equation}{section}
\NewTheorem{Assumption}
\NewTheorem{Proposition}
\NewTheorem{Theorem}
\NewTheorem{Corollary}
\NewTheorem{Conjecture}
\NewTheorem{Lemma}
\theoremstyle{definition}
\NewTheorem{Definition}
\theoremstyle{remark}
\NewTheorem{Remark}
\NewTheorem{Remarks}

\newaliascnt{Example}{equation}
\aliascntresetthe{Example}

\newenvironment{Example}%
  {\refstepcounter{Example}\trivlist
   \item[\hskip\labelsep\theequation.~\textbf{Example}\space]
   \ignorespaces
  }{\unskip\nobreak\hfil%
    \penalty50\hskip2em\hbox{}\nobreak\hfil$\Diamond$%
    \parfillskip=0pt\finalhyphendemerits=0\penalty-100\endtrivlist
}
\newaliascnt{Examples}{equation}
\aliascntresetthe{Examples}

  {\refstepcounter{Examples}\trivlist
   \item[\hskip\labelsep\theequation.~\textbf{Examples}\space]
   \ignorespaces #1\enumerate
   }{\endenumerate\unskip\nobreak\hfil%
    \penalty50\hskip2em\hbox{}\nobreak\hfil$\Diamond$%
    \parfillskip=0pt\finalhyphendemerits=0\penalty-100\endtrivlist
}

{\catcode`\|=\active
  \gdef\set#1{\mathinner{\lbrace\,{\mathcode`\|"8000%
                                   \let|\midvert #1}\,\rbrace}}
}
\def\midvert{\egroup\,\mid\,\bgroup}

\def\({\big(}
\def\){\big)}

\let\<=\langle
\let\>=\rangle

\usepackage{tikz}
\usetikzlibrary{matrix}

\newcount\tableauRow\newcount\tableauCol
\newcommand\Tableau[2][-2]{
  \begin{tikzpicture}[scale=0.3,draw/.append style={thick,black},baseline=#1mm]
    \tableauRow=0
    \foreach \Row in {#2} {
       \tableauCol=1
       \foreach\k in \Row {
          \draw(\the\tableauCol,\the\tableauRow)+(-.5,-.5)rectangle++(.5,.5);
          \draw(\the\tableauCol,\the\tableauRow)node{\k};
          \global\advance\tableauCol by 1
       }
       \global\advance\tableauRow by -1
    }
  \end{tikzpicture}
}

\def\bi{\mathbf{i}}

\newcommand{\K}{\mathcal{K}}
\renewcommand{\O}{\mathcal{O}}

\newcommand\RootMr[1][\pm]{\sqrt{\vrule height 2mm width 0pt\smash{M_r'}}}
\newcommand\RootLr[1][\pm]{\sqrt{\vrule height 2mm width 0pt\smash{L_r'}}}

\newcommand{\Sn}[1][n]{\mathfrak{S}_{#1}}
\newcommand{\An}[1][n]{\mathfrak{A}_{#1}}

\renewcommand{\H}{\mathscr{H}}

\newcommand{\Parts}[1][n]{\mathcal{P}_{#1}}

\newcommand{\C}{\mathbb{C}}
\newcommand{\Z}{\mathbb{Z}}

\def\Std{\mathop{\rm Std}\nolimits}

\newcommand{\la}{\lambda}

\newcommand\fo[1][\bi]{f^\O_{#1}}

\newcommand{\s}{{\tt s}}
\renewcommand{\t}{{\tt t}}

\renewcommand{\u}{{\tt u}}
\renewcommand{\v}{{\tt v}}
\newcommand\balpha{{\boldsymbol{\alpha}}}
\newcommand{\bla}{{\boldsymbol{\la}}}
\newcommand{\bmu}{{\boldsymbol{\mu}}}
\newcommand{\ei}{{\mathbf{i}}}
\newcommand{\ej}{{\mathbf{j}}}

\DeclareMathOperator{\res}{res}

\makeindex
\def\Email#1{\email{\href{mailto:#1}{#1}}}
\keywords{Alternating groups, alternating Hecke algebras,
Khovanov-Lauda-Rouquier algebras, representation theory}
\subjclass[2000]{20C08, 20D06, 20C30}
\begin{document}
\title[Representations of alternating cyclotomic Hecke algebras]{Semisimple representations of alternating cyclotomic Hecke algebras}
\author{Clinton Boys}
\address{School of Mathematics and Statistics F07, University of
Sydney\newline\hspace*{3mm} NSW 2006, Australia.}
\Email{clinton.boys@sydney.edu.au}


\begin{abstract}
We define alternating cyclotomic Hecke algebras in higher levels as subalgebras of cyclotomic Hecke algebras under an analogue of Goldman's hash involution. We compute the rank of these algebras and construct a full set of irreducible representations in the semisimple case, generalising Mitsuhashi's results \cite{Mitsuhashi}, \cite{Mitsuhashi2}. 
\end{abstract}

\maketitle

\newcommand\Macro[2][-]{%
\ifx#1e$\backslash$begin\{#2\}$\dots\backslash$end\{#2\}
\else $\backslash$#2\ifx#1b$\{\dots\}$\fi
\fi}

\section*{Introduction}
Cyclotomic Hecke algebras have a significant presence in the literature because of their deep connections to the representation theory of symmetric groups, quantum groups, and Lie theory. Recently, profound connections to the theory of Khovanov-Lauda-Rouquier (KLR) algebras have also been discovered. 

When one uncovers a result concerning symmetric groups and their representations, it is natural to ask what are the consequences for the alternating group, the index-2 subgroup of $\Sn$. In a similar way, Mitsuhashi introduced subalgebras of Iwahori-Hecke algebras, which function as $q$-analogues of alternating groups. Continuining in this way leads naturally to a family of algebras we call alternating cyclotomic Hecke algebras, a subfamily of which were studied by Ratliff \cite{Ratliff}. 

In this paper we extend this natural line of investigation and study the semisimple representation theory of alternating cyclotomic Hecke algebras in detail. We start by considering cyclotomic Hecke algebras, which grew via generalisation from the symmetric and hyperoctahedral groups, as well as their Iwahori-Hecke algebras. We also discuss the seminormal form theory of these algebras, which dates back to Young \cite{Young} in the symmetric group case, and is due to Hoefsmit \cite{Hoefsmit} and Ariki-Koike \cite{AK} in general for cyclotomic Hecke algebras -- we follow the more recent approach of Hu and Mathas \cite{HM3}. We then define alternating cyclotomic Hecke algebras for arbitrary level as fixed-point subalgebras of cyclotomic Hecke algebras under the hash involution. There is a clear analogy here to symmetric group algebras and the sign map. 
In the final section we give a classification of irreducible representations for semisimple alternating cyclotomic Hecke algebras, including a dimension formula. This gives a significant improvement on the existing proof for levels one \cite{Mitsuhashi} and two \cite{Mitsuhashi2} by Mitsuhashi, while extending it to arbitrary level. 

The family of algebras defined in this paper are further studied in \cite{BM}, where a graded isomorphism theorem is proved, providing a firm link between the alternating cyclotomic Hecke algebras and the world of quiver Hecke algebras in the spirit of Brundan and Kleshchev \cite{BK:GradedKL}. Specifically, it is shown in \cite{BM} that for $\Lambda=\Lambda_0$, the alternating cyclotomic Hecke algebra is isomorphic to a homogeneous subalgebra of the cyclotomic quiver Hecke algebra. In particular this gives a $\Z$-grading on the modular group algebras of alternating groups. 


{\bf Acknowledgements.} The author was supported by an Australian Postgraduate Award. 
The author would like to thank his supervisor Andrew Mathas for his guidance throughout his PhD candidature.  A slightly modified version of this paper appears 
as Chapter 3 of the author's PhD thesis \cite{BoysThesis} at the University of Sydney.

\section{Notation and definitions}




\newcommand{\Zo}{\mathcal{Z}}

In this chapter we give the definition of cyclotomic Hecke algebras as abstract algebras with a presentation by generators and relations. These are the objects we will be considering subalgebras of in this paper. In order to do this we need some notation. Let $\Zo$ be a unital integral domain.  
  If $k\in\Z$ and $\xi\in\Zo^\times$, define
  the {\em quantum integer}
  $[k]_\xi$ by
  $$[k]_\xi=\begin{cases}
    \phantom{-(}(1+\xi+\dots+\xi^{k-1}),& \text{if }k\ge0,\\
          -(\xi^{-1}+\xi^{-2}+\dots+\xi^k),&\text{if }k<0.
        \end{cases}$$
We write $[k]$ for $[k]_\xi$ when there is no confusion over the quantum parameter; by the geometric series formula,
  $[k]=\frac{\xi^k-1}{\xi-1}$ provided $\xi\neq 1$ (if $\xi=1$, $[k]=k$).


We can now define cyclotomic Hecke algebras. Our definition follows Ariki-Koike,  although our algebras comprise a subfamily of those considered in \cite{AK}. 

\begin{Definition}[Cyclotomic Hecke algebras \cite{AK}]\label{cHadef}
Let $\Zo$ be a unital integral domain, let $\xi\in \Zo^\times$ and let $\boldsymbol{\kappa}=(\kappa_1,\kappa_2,\ldots,\kappa_\ell)$ be an $\ell$-tuple of  elements in $\Z$. 
The {\em cyclotomic Hecke algebra} $\H_n=\H_{n,\ell}(\Zo,\xi,\boldsymbol{\kappa})$ is the unital associative $\Zo$-algebra generated by $L_1,L_2,\ldots,L_n,T_1,T_2,\ldots,T_{n-1}$ subject to the relations
\begin{equation}\label{cycheckrelns}
\begin{aligned}
\prod_{i=1}^\ell (L_1-[\kappa_i]_\xi)&= 0\\
(T_r+1)(T_r-\xi)&= 0\\
L_rL_s&= L_sL_r\\
T_rT_s&= T_sT_r,\quad\mbox{if $\left|r-s\right|>1$}\\
T_rT_{r+1}T_r&= T_{r+1}T_rT_{r+1}\\
T_rL_s&= L_sT_r,\quad\mbox{if $s\neq r,r+1$}\\
L_{r+1}(T_r-\xi+1)&= T_rL_r+1.
\end{aligned}
\end{equation}
The $\ell$-tuple $(\kappa_1,\ldots,\kappa_\ell)\in \Z^\ell$ is called the {\em multicharge} of $\H_n$ and is closely related to a {\em dominant weight} $\Lambda=\Lambda(\boldsymbol{\kappa})$ of a Kac-Moody algebra. The quantum integers $[\kappa_i]$ for $1\leq i\leq \ell$ are the {\em cyclotomic parameters}. 
The elements $L_1,L_2,\ldots,L_n$ are called {\em Jucys-Murphy elements}. $\ell$ is called the {\em level} of the algebra. 
\end{Definition}

\begin{Definition}\label{quantchar}
For a cyclotomic Hecke algebra $\H_n$, let $e$ be the {\em quantum characteristic} of $\xi$ in $\Zo$, that is, 
\[
e=\min\{k>0\mid 1+\xi+\xi^2+\ldots+\xi^{n-1}=0\},
\]
or $e=\infty$ if no such $k$ exists. Note that if $\xi=1$ and $\Zo=F$ is a field of positive characteristic, $e=\mathrm{char}(F)$. If $\xi\neq 1$ and $\xi$ is a root of unity then $e$ is the multiplicative order of $\xi$. 
\end{Definition}

The following basis theorem of Ariki-Koike gives the first indication that cyclotomic Hecke algebras are interesting objects of study: their rank depends only on $\ell$ and $n$ (i.e. it is independent of the choices of parameters). Most importantly, this rank does not depend on the choice of $\xi$ or its quantum characteristic. 

\begin{Theorem}[Ariki-Koike \cite{AK}]\label{AKbasis}
The cyclotomic Hecke algebra $\mathscr{H}_{n,\ell}(\Zo,\xi,\boldsymbol{\kappa})$ is a free $\Zo$-module with basis 
\[
\{L_1^{\gamma_1}L_2^{\gamma_2}\cdots L_n^{\gamma_n}T_\omega\mid 0\leq \gamma_i<\ell\; \mbox{and}\; \omega\in\Sn\},
\]
where $T_\omega=T_{i_1}T_{i_2}\cdots T_{i_k}$ if $s_{i_1}s_{i_2}\cdots s_{i_k}$ is any reduced expression for $\omega$. In particular its rank as a $\Zo$-module is 
\[\mathrm{rk}_\Zo(\mathscr{H}_{n,\ell}(\Zo,\xi,\boldsymbol\kappa))=\ell^nn!.\]
\end{Theorem}

In order to discuss the representation theory, semisimple and otherwise, of cyclotomic Hecke algebras, we need to introduce the combinatorial framework of (multi-) partitions and (multi-) tableaux. 

\begin{Definition}
A {\em partition} of $n\geq 0$ is a weakly decreasing sequence $\la=(\la_1,\la_2,\ldots)$ of non-negative integers which sum to $n$. An $\ell$-{\em multipartition} of $n$ is an $\ell$-tuple $\boldsymbol\la=(\la^{(1)},\la^{(2)},\ldots,\la^{(\ell)})$ of partitions whose total sum is $n$. When writing multipartitions, we omit trailing zeroes,  group repeated integers with exponents and separate components with bars. We write $\Parts$ for the set of partitions of $n$; $\la\vdash n$ means $\la\in\Parts$.  We write $\Parts^\ell$ for the set of $\ell$-multipartitions of $n$; $\boldsymbol\la\vdash_\ell n$ means $\boldsymbol\la\in\Parts^\ell$. 
\end{Definition}

%

\begin{Definition}
For a multipartition $\boldsymbol\la\in\Parts^\ell$, a $\boldsymbol\la$-{\em tableau} is a bijective filling of the boxes of $\boldsymbol\la$ with the numbers $1,2,\ldots,n$. A $\boldsymbol\la$-tableau is {\em standard} if the entries increase along rows and down columns within each constituent diagram. The collection of standard tableaux with $n$ boxes is written $\Std(\Parts^\ell)$; the collection of those of {\em shape} $\bla$ (i.e. those tableaux such that deleting all the numbers from the diagram recovers the diagram of the multipartition $\bla$) is written $\Std(\bla)$. 
\end{Definition}

\begin{Definition}\label{ordersdef}
We define the {\em dominance order} on $\Parts^\ell$ by writing $\bla\unrhd\bmu$, read as $\bla$ dominates $\bmu$, if
\[
\sum_{k=1}^{r-1}\left|\la^{(k)}\right|+\sum_{j=1}^i \la^{(r)}_j\geq \sum_{k=1}^{r-1}\left|\mu^{(k)}\right|+\sum_{j=1}^i \mu_j^{(r)}
\]
for all $1\leq r\leq \ell$ and $i\geq 1$. This is a partial order which gives $(\Parts^\ell,\unrhd)$ the structure of a poset. We can extend the dominance ordering to the set $\Std(\Parts^\ell)$ by defining $\s\unrhd\t$ if $\mathrm{sh}(\s\!\downarrow_m)\unrhd \mathrm{sh}(\t\!\downarrow_m)$ for all $1\leq m\leq n$, where by $\mathrm{sh}(\u)$ we mean the shape of the tableau $\u$, and by $\u\!\downarrow_k$ we mean the tableau with $k$ boxes obtained from $\u$ by deleting entries $k+1,k+2,\ldots,n$. 
\end{Definition}


There are two special tableaux for each multipartition $\bla\in\Parts^\ell$ which are used frequently: the {\em initial} $\bla$-tableau $\t^\bla$ which contains the entries $1,2,3,\ldots,n$ increasing along rows starting from $\la^{(1)}$, and the {\em final} $\bla$-tableau $\t_\bla$ which contains the same entries increasing down columns, starting from $\la^{(\ell)}$.

\begin{Definition}\label{conjdef}
For $\la\in\Parts$, the {\em conjugate} partition $\la'$ is the partition with
\[
\la'_j=\#\{i\geq 1\mid \la_i\geq j\}. 
\]
In terms of diagrams, $\la'$ is the diagram of $\la$ with rows and columns swapped. The {\em conjugate} $\bla'$ of a multipartition $\bla=(\la^{(1)},\ldots,\la^{(\ell)})\in\Parts^\ell$ is
\[
\bla'=(\la^{(\ell)}\;\!',\ldots,\la^{(1)}\;\!'). 
\]
For tableaux, conjugation $\t\mapsto\t'$ is defined by interchanging rows and columns; the {\em conjugate} of the multitableau $\t=(\t^{(1)},\ldots,\t^{(r)})$ is $\t'=(\t^{(\ell)}\;\!\!',\ldots,\t^{(1)}\;\!\!')$. Notice that
\begin{equation}\label{initialfinalconj}
(\t^\bla)'=\t_{\bla'}. 
\end{equation}
\end{Definition}

The following lemma gives the important relationship between the dominance ordering and the conjugation involution; a proof in the $\ell=1$ case is given in \cite[(1.11)]{MacDonald}; the general case follows from this. 

\begin{Lemma}\cite[(1.11)]{MacDonald}\label{conjinv}
Conjugation reverses the dominance order on multipartitions and multitableaux, that is, $\bla\unrhd\bmu$ if and only if $\bmu'\unrhd\bla'$ and $\s\unrhd\t$ if and only if $\t'\unrhd\s'$. 
\end{Lemma}

\begin{Definition}\label{resdef}
For $\t\in\Std(\bla)$, where $\bla\in\Parts^\ell$, the {\em content} of $k$ in $\t$, for $1\leq k\leq n$, is the number
\[
c_\t(k)=\kappa_l+c-r,
\]
where $k$ appears in component $l$ in column $c$ and row $r$ of $\t$; the {\em content sequence} of $\t$ is the $n$-tuple
\[
c(\t)=(c_\t(1),c_\t(2),\ldots,c_\t(n)). 
\]
The $e$-{\em residue} of $k$ in $\t$ is the residue of $c_\t(k)$ modulo $e$:
\[
\res_\t(k)=c_\t(k)+e\Z\in\Z/e\Z
\]
and the {\em residue sequence} of $\t$ is the $n$-tuple
\[
\ei_\t=(\res_\t(1),\res_\t(2),\ldots,\res_\t(n))\in (\Z/e\Z)^n. 
\]
Finally, for $\ei\in (\Z/e\Z)^n$, $\Std(\ei)$ is the set $\{\t\in\Std(\Parts^\ell)\mid \ei_\t=\ei\}$. 
\end{Definition}

Using the Jucys-Murphy elements and the notation from this section, we may produce a full set of mutually orthogonal idempotents for the Hecke algebras.  
For $\t\in\Std(\bla)$ for $\bla\in\Parts^\ell$, define
\begin{equation}\label{ssidem}
F_\t=\prod_{k=1}^n \prod_{\substack{\s\in\Std(\bla)\\ \res_\s(k)\neq \res_\t(k)}}\frac{L_k-[\res_s(k)]_\xi}{[\res_\t(k)]_\xi-[\res_\s(k)]_\xi}. 
\end{equation}

%
By \eqref{cycheckrelns}, the Jucys-Murphy elements commute so there is no ambiguity in the order of factors in \eqref{ssidem}. 
\begin{Proposition}[Ariki's semisimplicity criterion \cite{Ariki}]\label{arikithm}
Let $\H_n=\H_{n,\ell}(F,\xi,\boldsymbol\kappa)$ be a cyclotomic Hecke algebra with $e>2$, where $F$ is a field. Then 
$\H_n$ is a semisimple $F$-algebra if and only if the element
\begin{equation}\label{arikissc}
P_\H=P_\H(F,\xi,\boldsymbol\kappa)=[1]_\xi [2]_\xi\cdots [n]_\xi \prod_{1\leq r<s\leq \ell}\prod_{-n<d<n}[\kappa_r+d-\kappa_s]_\xi 
\end{equation}
is nonzero. 
Moreover, if $P_\H$ is nonzero then the collection $\{F_\t\mid\t\in\Std(\Parts^\ell)\}$ is a complete set of pairwise orthogonal idempotents for $\H_n$. 
\end{Proposition}

\section{The seminormal form}

Young \cite{Young} introduced the seminormal form for symmetric group algebras to give a particularly elegant description of their ordinary representation theory (i.e. over a field of characteristic zero). We will see the seminormal form arise as a basis of simultaneous eigenvectors for the Jucys-Murphy elements. 
As we will see, this basis, which is intimately linked to the semisimple representation theory of cyclotomic Hecke algebras, is particularly well-adapted to the computations we will perform in later chapters.


The following lemma states the well-known result that tableaux are uniquely determined by their content sequences (which are the same as their residue sequences for $e>n$) in the semisimple case.  

  \begin{Lemma}\cite[Lemma 3.34]{Mathas}\label{L:separation}
    Suppose that $\s,\t\in\Std(\Parts^\ell)$ and that $P_\H$ is nonzero. Then
    $\s=\t$ if and only if $[c_r(\s)]=[c_r(\t)]$,
    for $1\le r\le n$. In particular $\s=\t$ if and only if $\ei_\s=\ei_\t$. 
  \end{Lemma}

Recall that if $\t\in\Std(\Parts^\ell)$ and $1\le r\le n$ then $c_r(\t)\in\Z$ is the
content of~$r$ in~$\t$. 

\begin{Definition}
Define the integer $\rho_r(\t)\in\Z$ by
\begin{equation}\label{E:rhodef}
\rho_r(\t)=c_r(\t)-c_{r+1}(\t).
\end{equation}
$\rho_r(\t)$ is called the {\em axial distance} from $r+1$ to $r$ in~$\t$.
\end{Definition}

The following definition makes precise the freedom allowed when choosing the coefficients in Young's seminormal form. 

\begin{Definition}[\protect{Hu-Mathas~\cite[\S3]{HM3}}]
  \label{D:alphaSNCS}
  A {\em $*$-seminormal coefficient system} is a set of scalars
  $\balpha=\{\alpha_r(\s)\mid 1\leq r<n\text{ and } \s\in\Std(\Parts^\ell)\}$
  in~$\Zo$ such that 
  \begin{itemize}
  \item[(i)] for $\t\in\Std(\Parts^\ell)$ we have
  \begin{equation}\label{sncs1}
  \alpha_k(\t)\alpha_m(s_k\cdot\t)=\alpha_m(\t)\alpha_k(s_m\cdot \t)
  \end{equation}
for $1\leq k,m\leq n$ if $\left|k-m\right|>1$
\item[(ii)] for $\t\in\Std(\Parts^\ell)$ and $1\leq r\leq n-2$ we have 
  \begin{equation}\label{sncs2}
  \alpha_r(s_{r+1}s_r\t)\alpha_{r+1}(s_r\t)\alpha_r(\t)
      =\alpha_{r+1}(s_{r+1}s_r\t)\alpha_r(s_{r+1}\t)\alpha_{r+1}(\t)
  \end{equation}
 \item[(iii)] for $\t\in\Std(\Parts^\ell)$ and $1\le r<n$ and $\v=(r,r+1)\t$ then $\alpha_r(\s)=0$ if $\v\notin\Std(\Parts^\ell)$
  and otherwise
  \begin{equation}\label{E:sncf}
      \alpha_r(\t)\alpha_r(\v)
           =\frac{[1+\rho_r(\t)][1+\rho_r(\v)]}{[\rho_r(\t)][\rho_r(\v)]}.
  \end{equation}
  \end{itemize}
\end{Definition}

As noted in \cite[\S3]{HM3}, examples of seminormal
coefficient systems for the symmetric groups date back to Young in
1901~\cite{Young}. For example, we can take
$\alpha_r(\t)=\tfrac{[1+\rho_r(\t)]}{[\rho_r(\t)]}$, whenever
$\t, s_r\cdot\t\in\Std(\Parts^\ell)$. This is the choice made by James in \cite{James}. We give another example now, to motivate the particular choice we will make in later chapters. 



\begin{Example}\label{Ex:AltCS}
  Suppose that $\Zo=\K$ is a field which contains $\xi\in\K^\times$, $\sqrt{\xi}$ and square roots
  $\sqrt{-1}$ and $\sqrt{[h]_\xi}$, for $1\le|h|\le n$, such that $\sqrt{[-h]_\xi}=\sqrt{-1}(\sqrt{\xi})^h\sqrt{[h]_\xi}$, for
  $1<h\le n$. Define
  $$\alpha_r(\t)=\begin{cases}
     \frac{\sqrt{-1}\sqrt{[1+\rho_r(\t)]}\sqrt{[1-\rho_r(\t)]}}%
     {[\rho_r(\t)]}&\text{if }s_r\cdot\t\in\Std(\Parts^\ell)\text{ and }c_r(\t)>0\\
\frac{-\sqrt{-1}\sqrt{[1+\rho_r(\t)]}\sqrt{[1-\rho_r(\t)]}}%
     {[\rho_r(\t)]}&\text{if }s_r\cdot\t\in\Std(\Parts^\ell)\text{ and }c_r(\t)<0   \end{cases}
  $$
  and $\alpha_r(\t)=0$ if $s_r\cdot\t$ is not standard. 
  One can easily check that these scalars satisfy the requirements of a seminormal coefficient system. 
\end{Example}

\begin{Definition}
Let $*$ be the unique involutive anti-automorphism of $\H_n$ which fixes the generators $L_1,T_1,\ldots,T_{n-1}$. 
\end{Definition}

\begin{Definition}\label{snbdef}
A basis $\{f_{\s\t}\mid\s,\t\in\Std(\bla)\;\mbox{for}\;\bla\in\Parts^\ell\}$ of $\H_n=\H_{n,\ell}(\Zo,\xi,\boldsymbol{\kappa})$ is a {\em seminormal
basis} if 
\[
L_kf_{\s\t}=[c_k(\s)]f_{\s\t}\quad\mbox{and}\quad f_{\s\t}L_k=[c_k(\t)]f_{\s\t},
\]
for all $\s,\t\in \Std(\bla)$ with $\bla\in\Parts^\ell$ and
$1\le k\le n$. The above basis is a {\em $*$-seminormal basis} if in addition
$f_{\s\t}^*=f_{\t\s}$, for all $\s,\t\in\Std(\bla)$ with $\bla\in\Parts^\ell$. 
\end{Definition}

We can use seminormal coefficient systems to provide a complete description of the semisimple representation theory of cyclotomic Hecke algebras. Recall that the element $P_\H$ from Proposition \ref{arikithm} gives us a criterion for semisimplicity. 

\begin{Theorem}[\protect{%
  The Seminormal Form~\cite[Theorem 3.14]{HM3}}]
  \label{T:SeminormalForm}
  Suppose that $\K$ is a field in which $P_\H$ is nonzero and which contains a seminormal coefficient
  system $\balpha$  for $\H_n(\K)$. Then 
  \begin{itemize}
  \item[(i)] $\H_n(\K)$ has a $*$-seminormal basis
  $\{f_{\s\t}\mid\s,\t\in\Std(\bla)\;\mbox{for}\;\bla\in\Parts^\ell\}$, such that
  \[
  f_{\s\t}^*=f_{\t\s},\quad L_kf_{\s\t}=[c_k(\s)]f_{\s\t}\quad\text{and}\quad
    T_rf_{\s\t}=\alpha_r(\s)f_{\u\t}-\frac{1}{[\rho_r(\s)]}f_{\s\t},
  \]
  where $\u=s_r\cdot\s$ (and $f_{\u\t}=0$ if $\u$ is not standard);
  \item[(ii)] there exist non-zero scalars $\gamma_\t\in\K$, for $\t\in\Std(\Parts^\ell)$ such that
  \[
   f_{\s\t}f_{\u\v}=\delta_{\t\u}\gamma_\t f_{\s\v};
  \]
  \item[(iii)]  $\{\tfrac1{\gamma_\t}f_{\t\t}\mid\t\in\Std(\Parts^\ell)\}$ is a
  complete set of pairwise orthogonal primitive idempotents for $\H_n(\K)$, and
  \item[(iv)] the $*$-seminormal basis $\{f_{\s\t}\mid\s,\t\in\Std(\bla)\;\mbox{for}\;\bla\in\Parts^\ell\}$ is
  uniquely determined by the $*$-seminormal coefficient system $\balpha$ together with
  the scalars $\{\gamma_{\t^\bla}\mid\bla\in\Parts^\ell\}$.
  \end{itemize}
\end{Theorem}



Finally, we will need the following identity relating the $\alpha$ and $\gamma$ coefficients. 

\begin{Corollary}{\cite[Corollary 3.17]{HM3}}
  \label{L:alphagamma}
Suppose that $\t\in \Std(\Parts^\ell)$ is such that $\u=s_r\cdot\t\in\Std(\Parts^\ell)$, where $1\le r<n$. Then
$\alpha_r(\u)\gamma_\t=\alpha_r(\t)\gamma_\u$.
\end{Corollary}

By Theorem \ref{T:SeminormalForm}, if $\t\in\Std(\Parts^\ell)$ then
$F_\t=\tfrac1{\gamma_\t}f_{\t\t}$ is a primitive idempotent in $\H_n$. As we saw in \eqref{ssidem}, there is an explicit formula for this idempotent 
which in particular is independent of the choice of seminormal basis. Again by Theorem \ref{T:SeminormalForm}, as an $(\mathscr{L},\mathscr{L})$-bimodule $\H_n$ decomposes as
\begin{equation}\label{E:LLBimodule}
  \H_n=\bigoplus_{\substack{\bla\in\Parts^\ell\\\s,\t\in\Std(\bla)}}H_{\s\t},
\end{equation}
where $H_{\s\t}=\K f_{\s\t}$, and where $\mathscr{L}$ is the commutative subalgebra of $\H_n$ generated by the Jucys-Murphy elements. Equivalently,
$$H_{\s\t}=\{h\in\H_n\mid L_k h=[c_k(\s)]h\text{ and }
       hL_k=[c_k(\t)]h\text{ for }1\le k\le n\},$$
for $\s,\t\in\Std(\bla)\;\mbox{with}\;\bla\in\Parts^\ell$.



The following easy corollary of Theorem \ref{T:SeminormalForm} is the mechanism by which we perform many of the computations in this paper. It allows us to prove identities in cyclotomic Hecke algebras by comparing coefficients of their actions on the seminormal basis. Importantly, we can also use it prove identities in the {\em non-semisimple} versions of thse algebras.

\begin{Corollary}\label{identdecomp}
In the cyclotomic Hecke algebra $\H_n(\K)$, the identity element $1\in\H_n$ decomposes as 
\[
1=\sum_{\t\in\Std(\Parts^\ell)}\frac{1}{\gamma_\t}f_{\t\t}.
\]
\end{Corollary}

\section{Cyclotomic Hecke algebras with symmetric multicharges}

In this section we define the particular subfamily of cyclotomic Hecke algebras in which we are interested. This subfamily allows for the definition of an involution whose fixed-point subalgebra is our main topic. This construction will generalise the construction of the alternating group algebra as a fixed-point subalgebra of the symmetric group algebra.  

Recall the Jucys-Murphy elements from \eqref{cycheckrelns}. For $k=1,2,\ldots,n$, let
\[
\widetilde{L_k}=\xi^{1-k}T_{k-1}T_{k-2}\cdots T_2T_1L_1T_1T_2\cdots T_{k-2}T_{k-1}.
\]
These may be referred to as the {\em affine} Jucys-Murphy elements in $\H_n$. 
If $\xi\neq 1$, it follows by induction on $k$ that the two different definitions are related by
\begin{equation}\label{E:affineJM}
{L_k}=\frac{\widetilde{L}_k-1}{\xi-1}.
\end{equation}

\begin{Definition}
For a multicharge $\boldsymbol\kappa=(\kappa_1,\kappa_2,\ldots,\kappa_\ell)$, the multicharge $\boldsymbol\kappa'$ is defined by $(-\kappa_\ell,-\kappa_{\ell-1},\ldots,-\kappa_1)$. 
\end{Definition}

The involution we are interested in is the hash map, originally defined by Goldman, which we now define. Unlike the sign map for symmetric group algebras, it is not immediately obvious that the map we want is a homomorphism. 

\begin{Proposition}\label{hashprop}
Let $\Zo$ be a unital integral domain and let $\H_{n,\ell}(\Zo,\xi,\boldsymbol{\kappa})$ be a cyclotomic Hecke algebra. Then there is a unique algebra homomorphism 
\[
\#:\mathscr{H}_{n,\ell}(\Zo,\xi,\boldsymbol{\kappa})\to\mathscr{H}_{n,\ell}(\Zo,\xi,\boldsymbol{\kappa}')
\]
satisfying
\[
T_i\mapsto -\xi T_i^{-1}\quad\mbox{for $i=1,2,\ldots,n-1$.}
\]
and
\[
\begin{array}{ll}L_1\mapsto -L_1,&\mbox{if $\xi=1$}\\
\widetilde{L_1}\mapsto \widetilde{L_1}^{-1},&\mbox{if $\xi\neq 1$}\end{array}
\]
\end{Proposition}
\begin{proof}
We prove the result for the algebra $\widetilde{\H_n}=\H_{n,\ell}(\overline{\Z},\xi,\boldsymbol\kappa)$, where $\overline{\Z}=\Z[\xi,\xi^{-1},\boldsymbol\kappa]$; it is clear by Definition \ref{cHadef} that we can then base-change to our arbitrary integral domain $\O$ by $\H_n(\O)\cong \H_n(\overline{\Z})\otimes_{\overline{\Z}}\O$ to obtain the result in general. 

Note that $\widetilde{L_1},T_1,T_2,\ldots,T_{n-1}$ is also a generating set for the algebra $\widetilde{\H_n}$ when $\xi\neq 1$. Indeed, Ariki-Koike give such a presentation for $\H_n$ in \cite[Definition 3.1]{AK} which, for now ignoring the first relation, has all the same relations as in Definition \ref{cycheckrelns} but with the final relation replaced with the relation
\[
T_1\widetilde{L_1}T_1\widetilde{L_1}=\widetilde{L_1}T_1\widetilde{L_1}T_1.
\] 
Using this presentation, the fact that $\#$ preserves all relations in \eqref{cycheckrelns} except the $L_1$ eigenvalue relation are quick checks that we leave to the reader. For the remaining relation we split into two cases. First suppose $\xi=1$. Then, since $L_0=T_1$, 
\[
\Bigl(\prod_{i=1}^\ell (L_1-\kappa_i)\Bigr)^\#= \prod_{i=1}^\ell(-L_1-\kappa_i)= 0
\]
since if $\kappa_i$ appears in $\boldsymbol\kappa=\boldsymbol\kappa'$, so does $-\kappa_i$. On the other hand, if $\xi\neq 1$, let $\mathcal{K}$ be the field of fractions of $\overline{\Z}_{(\xi)}$ and let $\{f_{\s\t}\mid \s,\t\in\Std(\bla)\mbox{ for }\bla\in\Parts^\ell\}$ be a seminormal basis for the semisimple $\K$-algebra $\widetilde{\H_n}(\mathcal{K})$. Then $\widetilde{L_1}=(\xi-1){L_1}+1$ and so for a standard tableau $\t$, $\widetilde{L_1}f_{\t\t}=\xi^{c_1(\t)}f_{\t\t}=\xi^{\kappa_i}f_{\t\t}$ if 1 appears in component $i$ of $\t$. So, working in $\widetilde{\H_n}(\mathcal{K})$,
\begin{align*}
\Bigl(\prod_{i=1}^\ell (L_1-[\kappa_i])\Bigr)^\#f_{\t\t}&= \Bigl[\prod_{i=1}^\ell \Bigl(\frac{\widetilde{L_1}-1}{\xi-1}-[\kappa_i]\Bigr)\Bigr]^\#f_{\t\t}\\
&= \frac{1}{(\xi-1)^\ell}\Bigl(\prod_{i=1}^\ell(\widetilde{L_1}-\xi^{\kappa_i})\Bigr)^\#f_{\t\t}\\
&= \frac{1}{(\xi-1)^\ell}\prod_{i=1}^\ell(\widetilde{L_1}^{-1}-\xi^{\kappa_i})f_{\t\t}\\
&=0
\end{align*}
since $\boldsymbol\kappa=\boldsymbol\kappa'$ and $\widetilde{L_1}^{-1}f_{\t\t}=\xi^{-\kappa_i}f_{\t\t}$ if 1 appears in component $i$ of $\t$. The result now follows from Corollary \ref{identdecomp} since ${\H_n}(\overline{\Z})\hookrightarrow \widetilde{\H_n}(\mathcal{K})$. 
%
\end{proof}

\begin{Definition}\label{hashdef}
The map $\#:\mathscr{H}_{n,\ell}(\Zo,\xi,\boldsymbol{\kappa})\to\mathscr{H}_{n,\ell}(\Zo,\xi,\boldsymbol{\kappa}')$ from Proposition \ref{hashprop} is called the {\em hash map}. 
\end{Definition}

\begin{Remark}
The hash map was originally defined by Goldman \cite[Theorem~5.4]{Iwahori:Hecke} in level one; we have extended the definition to higher levels. 
\end{Remark}



In the same way that the alternating group algebra may be considered as the subalgebra of the symmetric group algebra of point fixed by the sign map, our goal now is to use the hash map we have just defined to give an analogous construction for cyclotomic Hecke algebras with symmetric parameters. In order to proceed, we will need to perform a number of calculations, computing the images of various elements we have defined so far under the hash involution from Definition \ref{hashdef}.

\begin{Lemma}
For $1\leq k\leq n$ we have $\widetilde{L_k}^\#=\widetilde{L_k}^{-1}$. 
\end{Lemma}
\begin{proof}
Using the final relation from \eqref{cycheckrelns} and noting $\xi T_r^{-1}=T_r-\xi+1$ we compute that $L_{r+1}=\xi^{-1}T_rL_rT_r+\xi^{-1}T_r$ and so we obtain an inductive formula for $\widetilde{L_{r+1}}$ as
\[
\widetilde{L_{r+1}}=\xi^{-1}T_r\widetilde{L_r}T_r. 
\]
Hence
\begin{align*}
\widetilde{L_k}^\#&=(\xi^{-1}T_{k-1}\widetilde{L_{k-1}}T_{k-1})^\#\\
&= \xi^{-1}\xi^2T_{k-1}^{-1}\widetilde{L_{k-1}}^{-1}T_{k-1}^{-1}\\
&= \widetilde{L_k}^{-1}
\end{align*}
by induction. 
\end{proof}

For the remainder of this section we need to be more careful with our choices of rings and fields. The following definition gives us the freedom to use the seminormal form to make meaningful statements for our algebras in general. 

\begin{Definition}\cite[Definition 4.1]{HM3}\label{eidem}
Let $\mathcal{K}$ be a field in which $P_\H$ is nonzero. Let $\O$ be a subring of $\mathcal{K}$ and $t\in\O^\times$. Then $(\O,t)$ is an {\em $e$-idempotent subring} of $\mathcal{K}$ if the following hold:
\begin{itemize}
\item[(i)] $[k]_t$ is invertible in $\O$ whenever $k\not\equiv 0$ (mod $e$) for $k\in\Z$; and 
\item[(ii)] $[k]_t\in \mathcal{J}(\O)$ whenever $k\in e\Z$,
\end{itemize}
where $\mathcal{J}(\O)$ is the Jacobson radical of $\O$, i.e. the intersection of all its maximal ideals. 
\end{Definition}

\label{fields}



Let $F$ be an arbitrary field with $\xi\in F^\times$ such that the quantum characteristic of $\xi$ in $F$ is $e>2$, and let $(\O,t)$ be an $e$-idempotent subring of a field $\mathcal{K}$ such that
\begin{itemize}
\item[(i)] $\K$ contains a seminormal coefficient system for $\H_n$; and
\item[(ii)] $F=\O/\mathfrak{m}$ for some maximimal ideal $\mathfrak{m}$ of $\O$ and $\xi=t+\mathfrak{m}$.
\end{itemize}

\begin{Remark}
It is shown in \cite[Example 4.2]{HM} that $e$-idempotent subrings exist. 
Note that by (ii) above, $\H_n(F)\cong \H_n(\O)\otimes_\O F$ and, since $\O$ is a subring of $\K$, $\H_n(\K)\cong \H_n(\O)\otimes_\O \K$. 
\end{Remark}

\begin{framed}
For the rest of this section, let us fix the notation above, together with a seminormal basis $\{f_{\s\t}\mid \s,\t\in\Std(\bla)\mbox{ for }\bla\in\Parts^\ell\}$ for $\H_n(\mathcal{K})$ and a seminormal coefficient system $\boldsymbol\alpha$. In particular, for the remainder of this section, $e>2$. 
\end{framed}

We now perform a number of calculations in the semisimple cyclotomic Hecke algebra $\H_n(\mathcal{K})$. It is clear from Proposition \ref{hashprop} that if $\H_n$ is a cyclotomic Hecke algebra with multicharge $\boldsymbol\kappa$ such that $\boldsymbol\kappa=\boldsymbol\kappa'$, then $\#$ is an {\em involution} of $\H_n$. We use this assumption implicitly in many calculations below. 

\begin{Lemma}\label{L:LFHash}
  Suppose that $\boldsymbol\kappa=\boldsymbol\kappa'$, that $1\le k\le n$ and $\s\in\Std(\Parts^\ell)$. Then
  $$L_k^\# f_{\s\s}=[c_k(\s')]f_{\s\s}.$$
\end{Lemma}
\begin{proof}
We give a proof if $\xi\neq 1$; the proof when $\xi=1$ is easier and we leave it to the reader.  By Theorem \ref{T:SeminormalForm}, $L_k f_{\s\t}=[c_k(\s)]f_{\s\s}$, so
  $\widetilde L_kf_{\s\s}=\xi^{c_k(\s)}f_{\s\s}$ by \eqref{E:affineJM}.
  By Definition \ref{hashdef} then,
  $\widetilde{L}_k^\#=\widetilde{L}_k^{-1}$. Therefore,
  \begin{align*}
  L_k^\# f_{\s\s}&=\frac{\tilde L_k^\#-1}{\xi-1}f_{\s\s}\\
                   &=\frac{\tilde L_k^{-1}-1}{\xi-1}f_{\s\s}\\
                   &=\frac{\xi^{-c_k(\s)}-1}{\xi-1}f_{\s\s}\\
                   &=[c_k(\s')]f_{\s\s},
                   \end{align*}
   where the last equality follows because $c_k(\s')=-c_k(\s)$ by Definitions \ref{conjdef} and \ref{resdef}. 
\end{proof}

\begin{Lemma}\label{L:Fhash}
Suppose that $\boldsymbol\kappa=\boldsymbol\kappa'$ and $\s\in \Std(\Parts^\ell)$. Then $F_\s^\#=F_{\s'}$.
\end{Lemma}
\begin{proof}
Since $F_\s=\frac{1}{\gamma_\s}f_{\s\s}$, applying Lemma \ref{L:LFHash} gives
\[
L_kF_\s^\#=(L_k^\#F_\s)^\#=([c_k(\s')F_\s)^\#=[c_k(\s')]F_\s^\#
\]
since the $\#$ map is an automorphism. Similarly $F_\s^\#L_k=[c_k(\s')]L_k$. Therefore, $F_\s^\#\in H_{\s'\s'}$ in the
decomposition of \eqref{E:LLBimodule}. As $F_\s$ is an idempotent, and $\#$ is
an algebra isomorphism, it follows that $F_\s^\#=F_{\s'}$ since this is the
unique idempotent in $H_{\s'\s'}=\K F_{\s'}$.
\end{proof}

\begin{Corollary}\label{C:ftthash}
  Suppose that $\boldsymbol\kappa=\boldsymbol\kappa'$ and $\s\in\Std(\Parts^\ell)$. Then
  \[
  f_{\s\s}^\#=\frac{\gamma_\s}{\gamma_{\s'}}f_{\s'\s'}.
  \]
\end{Corollary}

\begin{proof} Using Theorem \ref{T:SeminormalForm} and Lemma \ref{L:Fhash},
  $f_{\s\s}^\#=\frac1{\gamma_\s}F_\s^\#=\frac1{\gamma_\s}F_{\s'}
  =\frac{\gamma_{\s'}}{\gamma_\s}f_{\s'\s'}$.
\end{proof}


\begin{Definition}[Cyclotomic Hecke algebra with symmetric multicharge]\label{cHacpdef}
We say the cyclotomic Hecke algebra $\H_{n,\ell}(\Zo,\xi,\boldsymbol\kappa)$ has {\em symmetric multicharge} if $\boldsymbol\kappa=\boldsymbol\kappa'$.
\end{Definition}

 The symmetry of the multicharge which allowed us to perform the calculations above also gives us the following useful combinatorial lemma. 

\begin{Lemma}\label{minusres}
For a cyclotomic Hecke algebra $\H_n$ with symmetric multicharge, $\res_k(\t)\equiv -\res_k(\t')$ mod $e$ for all $\t\in\Std(\bla)$ for $\bla\in \Parts^\ell$ and $1\leq k\leq n$. 
\end{Lemma}
\begin{proof}
Let  $\bla\in \Parts^r$. Then for $\t\in\Std(\bla)$ and $1\leq k\leq n$, we have
\begin{align*}
\res_k(\t)&=\kappa_\ell+r_k(\t)-c_k(\t)\\
&\equiv -\kappa_{r-\ell+1}-r_k(\t')+c_k(\t')\\
&=-(\kappa_{r-\ell+1}+r_k(\t')-c_k(\t'))\\
&=-\res_k(\t'). 
\end{align*}
as required. \end{proof}

\begin{Remark}$\;$

\begin{itemize}
\item[(i)] Although our algebras $\H$ do not themselves depend on the choice of multicharge (only on the residues modulo $e$), the choice of seminormal coefficient system, and therefore the algebra $\H^\O$, does depends on this choice. Our notation reflects this dependence on $\kappa$. 
\item[(ii)] The reader may notice that for the proof of Proposition \ref{hashprop}, we could have made the slightly weaker assumption that $\kappa_i\in \boldsymbol\kappa$ implies $-\kappa_i\in \boldsymbol\kappa$; it is for Lemma \ref{minusres} that we need the multicharge to be symmetric in this particular way. This leads to the slightly uncomfortable reality that $\boldsymbol\kappa=(0,1,2)$ is not a symmetric 3-multicharge, but $\boldsymbol\kappa=(1,0,2)$ is. 
\end{itemize}
\end{Remark}

We want to study the subalgebra of a cyclotomic Hecke algebra $\H_n$ with a symmetric multicharge consisting of elements fixed by the hash involution; the study of these algebras 
will occupy the rest of this paper. 

\begin{Definition}[Alternating cyclotomic Hecke algebras]\label{acHadef}
Let $\Zo$ be a unital integral domain such that $2$ is invertible in $\Zo$. Fix $\xi\neq -1$ in $\Zo$. The {\em alternating cyclotomic Hecke algebra} $\H_n^\#=\H_{n,\ell}(\Zo,\xi,\boldsymbol{\kappa})^\#$ of type $(\ell,n)$ with Hecke parameter $\xi$ and symmetric cyclotomic parameters $\boldsymbol\kappa=\boldsymbol\kappa'$ is the fixed-point subalgebra of $\H_{n,\ell}(\Zo,\xi,\boldsymbol\kappa)$ under the $\#$ involution. 
\end{Definition}


\begin{Remark}$\;$

\begin{itemize}
\item[(i)] Mitsuhashi first studied alternating cyclotomic Hecke algebras in \cite{Mitsuhashi}, who called them {\em alternating Iwahori-Hecke algebras} by analogy with existing terminology for the $\ell=1$ case. For $q\in\C^\times$ with $q\neq -1$, Mitsuhashi \cite[Definition 4.1]{Mitsuhashi} defines a subalgebra of the Iwahori-Hecke algebra $\H_q(\Sn)=\H_{n,1}(\C,q,0)$, denoted by $\H_q(\An)$, which satisfies $\H_1(\An)\cong \C\An$. His definition is not by the hash involution, but can be shown to be equivalent to ours. A description of its representation theory is obtained for generic $q$ \cite[Theorem 5.5]{Mitsuhashi}, as well as a presentation by generators and relations which is a $q$-analogue of the presentation for $\C\An$ \cite[\S4]{Mitsuhashi}. 
\item[(ii)] Mitsuhashi also studied the case $\ell=2$ when $\boldsymbol\kappa=(1,e-1)$ \cite{Mitsuhashi2} (using a Clifford theory approach similar to \cite[\S1]{BM}). His algebra is isomorphic to our algebra $\H_{n,2}(\C,q,(1,-1))^\#$ and his results give explicit generators and relations for this family of alternating cyclotomic Hecke algebras \cite[Proposition 3.2]{Mitsuhashi2}. 
\end{itemize}
\end{Remark}

Our next goal is to determine the dimension of $\H_n^\#$ over a field. First we need some lemmas and some additional notation. For $\ei\in I^n$, define
\begin{equation}\label{idemdef}
f_\ei^\O=\sum_{\substack{\t\in\Std(\Parts^\ell)\\\res(\t)=\ei}} F_\t.
\end{equation}

These idempotents will be very important in later chapters, as they appear in a different guise when we study these algebras through a different looking glass. The following standard result shows that, despite ostensibly belonging to $\H_n(\mathcal{K})$, the idempotents $\fo$ actually belong to the $\O$-form of the algebra, justifying their notation. 

\begin{Lemma} \cite[Lemma 4.5]{HM3}, \cite{MurphyIdem}\label{eio}
For $\ei\in I^n$, $\fo\in \H_n(\O)$. 
\end{Lemma}

These elements generalise Proposition \ref{arikithm} in the following sense. 

\begin{Proposition}\cite{Grojnowski, MathasTilting} If the element $P_\H(F,\xi,\boldsymbol\kappa)$ is zero, then $\H_n(F)$ is a non-semisimple $F$-algebra and the collection $\{f_\ei^\O\otimes_\O1_F\mid \ei\in I^n\}$ is a complete set of pairwise orthogonal idempotents for $\H_n(F)$.   
\end{Proposition}

We now compute the image of the idempotent $\fo$ under the hash map. Given $\ei\in I^n$, define $-\ei\in I^n$ by
\begin{equation}\label{minusidef}
-\ei=(-i_1,-i_2,\ldots,-i_n).
\end{equation}

\begin{Lemma}\label{L:ehash}
  Suppose that $\bi\in I^n$. Then $(\fo)^\#=\fo[-\bi]$ in $\H_n(\O)$. 
\end{Lemma}

\begin{proof}
  First observe that $\s\in\Std(\bi)$ if and only if $\s'\in\Std(-\bi)$ by Lemma \ref{minusres}. Then by Lemma \ref{L:Fhash},
  \[
(\fo)^\#=\sum_{\s\in\Std(\bi)}F_\s^\#=\sum_{\s\in\Std(\bi)}F_{\s'}=\fo[-\bi]
  \]
  as claimed.
\end{proof}

We want to work with equivalence classes of $I^n$ under the involution on residue sequences defined in \eqref{minusidef}.  More precisely, for $\ei,\ej\in I^n$ let $\sim$ be the equivalence relation on $I^n$ generated by $\ei\sim \ej$ if $\ei=-\ej$. From this we obtain a partition of the set $I^n$ into equivalence classes of size 1 or 2; we denote the equivalence class containing a sequence $\ei$ by $[\ei]$, noting that, since $e>2$, 
\[
[\ei]=\left\{\begin{array}{ll}\{\ei\},&\mbox{if $\ei=\underbrace{(0,0,\ldots,0)}_{\text{$n$ zeroes}}$}\\ \{\ei,-\ei\},&\mbox{otherwise}\end{array}\right. 
\]
We denote the set of equivalence classes by $I^n_\sim$ and in each equivalence class we choose a representative $\ei^+\in[\ei]$. \label{eiplusdef}

We now give a dimension formula for alternating cyclotomic Hecke algebras.


\begin{Theorem}\label{ungradedaltdim}
Suppose that $F$ is a field of characteristic greater than two, and suppose that $\boldsymbol\kappa=\boldsymbol\kappa'$ is such that 
\[
\left|\{j\mid \kappa_j\equiv 0\;\mbox{mod $e$}\}\right|<n.
\]
Then the alternating cyclotomic Hecke algebra $\H_n(F)^\#$ has dimension $\dfrac{r^nn!}2$. 
\end{Theorem}
\begin{proof}
Choose a seminormal coefficient system in a field $\K$ with an $e$-idempotent subring $\O$ such that $\K\supset\O\to F$ as on p\pageref{fields}. Define an element
\[
\varepsilon=\sum_{\substack{[\ei]\in I^n_\sim\\\left|[\ei]\right|=2}} \Bigl(f_{\ei^+}^\O-f_{-\ei^+}^\O\Bigr)
\]
and note that, since $e>2$ and since the condition on $\boldsymbol\kappa$ disallows the sequence $(0,0,\ldots,0)$, $\left|[\ei]\right|=2$ for all $[\ei]\in I^n_\sim$. Then $\H_n^\#\cong \varepsilon\H_n^\#$ as $\O$-modules, since $\varepsilon^2=1$ and $\varepsilon^\#=-\varepsilon$ by Lemma \ref{L:ehash}. Writing $\fo[{[\ei]}]$ for $\fo[{\ei}^+]+\fo[{-\ei^+}]$, we see $\H_n\cong \H_n^\#\oplus \varepsilon\H_n^\#$, which gives the result by Theorem \ref{AKbasis} since any $x\in\H_n$ may be written as $x=\sum_{\ei\in I^n}x\fo=\sum_{[\ei]\in I^n_\sim}x\fo[{[\ei]}]$ and so we can write
\[
x=\frac12\sum_{[\ei]\in I^n_\sim}(x+x^\#)\fo[{[\ei]}]+ \frac12\varepsilon\sum_{[\ei]\in I^n_\sim}(x-x^\#)\fo[{[\ei]}]
\]
provided $\frac12\in\O$. The result now follows by tensoring with the field $F$. 
\end{proof}


\section{Semisimple representations of alternating cyclotomic Hecke algebras}

In this section we construct the semisimple representations of alternating cyclotomic Hecke algebras using the seminormal form from \S3.4. Our goal is to construct a full set of pairwise non-isomorphic irreducible modules for $\H_n(\mathcal{K})^\#$ for certain fields $\mathcal{K}$. This generalises results in \cite{Mitsuhashi} and \cite{Mitsuhashi2}. To begin with, we need to place some more restrictive conditions on our seminormal coefficient system $\boldsymbol\alpha$. 

\begin{Definition}\label{D:AltCS}
  An {\em alternating coefficient system} is a $*$-seminormal coefficient system
  $\balpha=\{\alpha_r(\s)\mid 1\leq r\leq n\;\mbox{ and }\s\in\Std(\Parts^r)\}$ such that $ \alpha_r(\s)=-\alpha_r(\s')$,
  for all $1\leq r<n$ and $\s\in\Std(\Parts^r)$.
\end{Definition}

\begin{Remark}
Example \ref{Ex:AltCS} shows that alternating seminormal coefficient systems exist.
\end{Remark}

\begin{framed}
For the remainder of this section, fix a field $\mathcal{K}$ with $\xi\in \K^\times$ with quantum characteristic $e>2$ and such that $P_\H$ is nonzero in $\mathcal{K}$, a symmetric multicharge $\boldsymbol\kappa$, a seminormal basis $\{f_{\s\t}\mid\s,\t\in\Std(\bla)\;\mbox{for}\;\bla\in\Parts^\ell\}$ for $\H_n(\mathcal{K})$ and an alternating seminormal coefficient system
$\balpha$.
\end{framed}
%


\begin{Remark}\label{alpha3}
Since we will need these coefficients explicitly later, we now compute $\alpha_r(\t)$ for some particularly important $r$ and $\t$. Note that by \eqref{E:sncf}, if $\t$ is any tableau with $2$ in the first row and $3$ in the first column, 
\[
\alpha_2(\t)\alpha_2(s_2\cdot\t)=\frac{[3][-1]}{[2][-2]}=\frac{-t^{-1}[3][1]}{-t^{-2}[2]^2}=\frac{t[3]}{[2]^2}.
\]
In order to satisfy the requirement of an alternating seminormal coefficient system, 
 we make the following choice for the sake of definiteness:
\begin{equation}\label{E:alphadef}
\alpha_2(\t)= \left\{\begin{array}{ll}\frac{\sqrt{-1}\sqrt{t}\sqrt{[3]}}{[2]},&\mbox{if 2 is in the first row of $\t$}\\ -\frac{\sqrt{-1}\sqrt{t}\sqrt{[3]}}{[2]},&\mbox{if 2 is in the second row of $\t$}\end{array}\right.
\end{equation}
This has the effect that $\alpha_2(\s)=\alpha_2(\widetilde{\s})$ for any tableaux $\s,\widetilde{\s}$ with $\s\downarrow_3=\widetilde{\s}\downarrow_3$.  
\end{Remark}


The following lemma continues some calculations we began in the previous section. 

\begin{Lemma}\label{L:futhash}
  Suppose that $\s,\u\in\Std(\bla)$ with $\bla\in\Parts^\ell$ are standard tableaux such that
  $\u=s_r\cdot\s$, where~$1\le r<n$. Then
\[
f_{\u\s}^\#=-\frac{\alpha_r(\s')\gamma_\s}{\alpha_r(\s)\gamma_{\s'}}f_{\u'\s'}.
\]
\end{Lemma}
\begin{proof}
By Theorem \ref{T:SeminormalForm},
$f_{\u\s}=\frac{1}{\alpha_r(\s)}\Big(T_r+\frac1{[\rho_r(\s)]}\Big)f_{\s\s}$.
Hence, using Definition \ref{hashdef} and Corollary \ref{C:ftthash} for the second
equality,
\begin{align*}
f_{\u\s}^\#&=\frac{1}{\alpha_r(\s)}\Big(T_r+\frac1{[\rho_r(\s)]}\Big)^\#f_{\s\s}^\#\\
&=\frac{\gamma_\s}{\alpha_r(\s)\gamma_{\s'}}\Big(-T_r+t-1+\frac1{[\rho_r(\s)]}\Big)f_{\s'\s'}\\
&=-\frac{\gamma_\s}{\alpha_r(\s)\gamma_{\s'}}\Big(T_r-\frac{t^{\rho_r(\s)}}{[\rho_r(\s)]}\Big)f_{\s'\s'}\\
&=-\frac{\gamma_\s}{\alpha_r(\s)\gamma_{\s'}}\Big(T_r+\frac1{[\rho_r(\s')]}\Big)f_{\s'\s'}
\end{align*}
since $[\rho_r(\s)]=-t^{\rho_r(\s)}[-\rho_r(\s)]=-t^{\rho_r(\s)}[\rho_r(\s')]$. Observe
that $\u'=s_r\s'$. Therefore, the result follows by another application of
Theorem \ref{T:SeminormalForm}.
\end{proof}

By Theorem \ref{T:SeminormalForm} any $*$-seminormal basis is uniquely determined by a
seminormal coefficient system and a choice of scalars
$\{\gamma_{\t^\bla}\mid\bla\in\Parts^\ell\}$. We now determine
these scalars for the seminormal basis
\[
\{f_{\s\t}^\#\mid\s,\t\in\Std(\bla)\;\mbox{for}\;\bla\in\Parts^\ell\}.
\] 


\begin{Proposition}\label{P:AltCS}
  The collection $\{f_{\s\t}^\#\mid\s,\t\in\Std(\bla)\;\mbox{for}\;\bla\in\Parts^\ell\}$ is the
  seminormal basis of $\H_n(\mathcal{K})$ determined by the seminormal coefficient system
  $$\{-\alpha_r(\s)\mid\s\in\Std(\bla),\; \bla\in \Parts^\ell\text{ and }1\le r<n\}$$
  together with the $\gamma$-coefficients $\{\gamma_{\t_\bla}\mid\bla\in\Parts^\ell\}$.
  That is, if $\s,\t\in\Std(\bla)$ for $\bla\in\Parts^\ell$ and $1\le r<n$ then
  $$T_r f_{\s\t}^\# = -\alpha_r(\s)f_{\u\t}^\#-\frac1{[\rho_r(\s')]}f_{\s\t}^\#,$$
  where $\u=s_r\cdot\s$. Moreover,
  $f_{\s\t}^\#f_{\u\v}^\#=\delta_{\t\u}\gamma_\t f_{\s\v}^\#$,
  for $\s,\t,\u,\v\in\Std(\bla)$ for $\bla\in \Parts^\ell$. 
\end{Proposition}

\begin{proof} Using
  Theorem \ref{T:SeminormalForm}, if $\s,\t\in\Std(\bla)$ for $\bla\in \Parts^\ell$ we compute
  \begin{align*}
    T_rf_{\s\t}^\# &= \Big(T_r^\# f_{\s\t}\Big)^\#
         = \Big((-T_r+t-1)f_{\s\t}\Big)^\#\\
        &=\Big(-\alpha_r(\s)f_{\u\t} +(t-1+\frac1{[\rho_r(\s)]})f_{\s\t}\Big)^\#\\
        &=\Big(-\alpha_r(\s)f_{\u\t}-\frac1{[\rho_r(\s')]}f_{\s\t}\Big)^\#\\
        &=-\alpha_r(\s)f_{\u\t}^\# -\frac1{[\rho_r(\s')]}f_{\s\t}^\#.
\end{align*}
Similarly,
$f_{\s\t}^\#f_{\u\v}^\#=(f_{\s\t}f_{\u\v})^\#=\delta_{\t\u}\gamma_\t f_{\s\v}^\#$.
By Theorem \ref{L:LFHash} $f_{\s\t}^\#\in H_{\s'\t'}$, so the $\alpha$-coefficient
corresponding to~$f_{\s\t}^\#$ is naturally indexed by~$\s'$ (and not by~$\s$).
Similarly, the labelling for the $\gamma$-coefficients involves conjugation
because $F_\t=\frac1{\gamma_{\t'}}f_{\t'\t'}^\#$ by Corollary \ref{C:ftthash}. Hence,
the result follows by Theorem \ref{T:SeminormalForm}.
\end{proof}


We now define Specht modules for semisimple cyclotomic Hecke algebras.

\begin{Definition}
Let $\bla\in \Parts^\ell$. The Specht module ${S^\bla}$ for the algebra $\H_n(\K)$ is the vector space with basis $\{f_\t\mid\t\in\Std(\bla)\}$ and with $\H_n(\K)$-action given by
\begin{align*}
L_kf_\t&= [c_k(\t)]f_\t\\
T_rf_\t&= \alpha_r(\t)f_\u-\frac{1}{[\rho_r(\t)]}f_\t
\end{align*}
for $1\leq k\leq n$ and $1\leq r<n$, where $\u=s_r\cdot\s$ (and $\alpha_r(\t)=0$ if $\u$ is not standard). 
\end{Definition}


\begin{Example}\label{spechex1}
Let us calculate the actions of the generators $T_1$ and $T_2$ on the Specht module ${S^{(21)}}$ for the algebra $\H_{3,1}(\C,1,0)\cong\C\mathfrak{S}_3$, noting that $\C$ clearly contains the required alternating seminormal coefficient system computed in Remark \ref{alpha3}. We see that $T_1$ and $T_2$ respectively act as the matrices
\[
\begin{pmatrix}1&0\\0&-1\end{pmatrix},\qquad \begin{pmatrix}\tfrac12&\tfrac{\sqrt{3}i}{2}\\ -\tfrac{\sqrt{3}i}{2}&-\frac12\end{pmatrix}
\]
which the reader can check square to the identity (the first column of each matrix corresponds to the vector $f_{12/3}$ and the second to $f_{13/2}$). 
\end{Example}


\begin{Theorem}\cite{AK}\label{ssspecht}
Let $\K$ be a field containing a seminormal coefficient system and such $P_\H$ is nonzero. Then for each $\bla\in\Parts^\ell$, ${S^\bla}$ is an irreducible $\H_n(\K)$-module. Moreover, the collection $\{{S^\bla}\mid\bla\in\Parts^\ell\}$ is a complete list of irreducible modules for the semisimple algebra $\H_n(\K)$. 
\end{Theorem}


\begin{Definition}
Let $\sim$ be the equivalence relation on $\Parts^\ell$ generated by $\bla\sim\bmu$ if $\bmu=\bla'$. We write $[\bla]$ for the equivalence class of $\bla\in\Parts^\ell$ under this equivalence relation, and $(\Parts^\ell)_\sim$ for the set of all equivalence classes. 
\end{Definition}

If $A$ is an $\O$-algebra and $B$ an $\O$-subalgebra of $A$, for an $A$-module $M$ we write $M\!\downarrow^A_B$ for the $B$-module obtained by restriction. 

\begin{Proposition}\label{ss1}
Let $\K$ be a field containing a seminormal coefficient system and such that $P_\H$ is nonzero. Let $\bla\in\Parts^\ell$ be such that $\bla\neq \bla'$. Then 
\[
{S^\bla}\!\downarrow^{\H_n(\K)}_{\H_n(\K)^\#}\cong {S^{\bla'}}\!\downarrow^{\H_n(\K)}_{\H_n(\K)^\#}
\]
as $\H_n(\K)^\#$-modules.  
\end{Proposition}
\begin{proof}
Define an map of vector spaces $\tau:S^\bla\to S_\bla$ by $v_\t\mapsto v_{\t'}$; this is clearly an isomorphism of vector spaces since $\left|\Std(\bla)\right|=\left|\Std(\bla')\right|$ by Lemma \ref{conjinv}. Then by Proposition \ref{P:AltCS}, $\tau(T_rv_\t)=T_r^\#v_{\t'}$ and by Lemma \ref{L:LFHash} $L_kv_\t=L_k^\#v_{\t'}$ so $\tau$ is an $\H_n(\K)^\#$-module isomorphism between the restricted modules. 
\end{proof}

\begin{Definition}
If $\bla\in\Parts^\ell$ is such that $\bla\neq\bla'$, we write $S^{[\bla]}$ for the $\H_n(\K)^\#$-module from Proposition \ref{ss1}. 
\end{Definition}

\begin{Corollary}\label{dim1}
Let $\K$ be a field containing a seminormal coefficient system and such that $P_\H$ is nonzero. Then for $\bla\in\Parts^\ell$ such that $\bla\neq \bla'$, $S^{[\bla]}$ is an irreducible $\H_n(\K)^\#$-module. 
\end{Corollary}
\begin{proof}
Let $v$ be a nonzero vector in $S^{[\bla]}= {S^\bla}\!\downarrow^{\H_n(\K)}_{\H_n(\K)^\#}$. Writing $v=\sum_{\t\in\Std(\bla)}r_\t v_\t$ as an arbitrary linear combination, we can surely choose some $\t$ with $r_\t\neq 0$. Then $r_\t f_\t=(F_\t+F_{\t'}) v\in \H_n^\# v$ and so $f_\t\in \H_n^\# v$. We now observe that it is possible to move from any one basis vector to another by applying a sequence of $\#$-invariant elements; one can check that for $\t\in\Std(\bla)$ with $\u=s_r\cdot\t$ standard, $(F_\u T_r+F_{\u'}T_r^\#)f_\t=\frac{\alpha_r(\t)}{\gamma_\u}f_\u$. We are done since we can clearly move all the way from the basis vector $f_{\t^\la}$ to the vector $f_{\t_\la}$ in such a way; the element we are acting by is $\#$-invariant by Lemma \ref{L:Fhash}. So $\H_n^\# v=S^{[\bla]}$ and we are done.
\end{proof}

\begin{Definition}
Let $\bla\in\Parts^\ell$ be such that $\bla=\bla'$ and suppose that 2 is invertible in $\K$. Denote by $S^\bla_+$ the vector space with basis $\{\frac12(f_\t+f_{\t'})\mid\t\in\Std(\bla)^+\}$ and by $S^\bla_-$ the vector space with basis $\{\frac12(f_\t-f_{\t'})\mid\t\in\Std(\bla)^+\}$.
\end{Definition}

Note that for the next proposition, it is important that we include the additional assumption that $\K$ contains an {\em alternating} seminormal coefficient system. 

\begin{Proposition}\label{dim2}
Let $\K$ be a field of characteristic greater than 2 containing an alternating seminormal coefficient system and such that $P_\H$ is nonzero. Then for $\bla\in\Parts^\ell$ with $\bla=\bla'$, $S^\bla_+$ and $S^\bla_-$ are irreducible $\H_n(\K)^\#$-modules. 
\end{Proposition}
\begin{proof}
Similarly to the proof of Corollary \ref{dim1}, 
suppose we need to move from $\frac12(f_\t\pm f_{\t'})$ to $\frac12(f_\u\pm f_{\u'})$, where $\u=s_r\cdot\t$ is standard. Then we compute, using Theorem \ref{T:SeminormalForm} and Proposition \ref{P:AltCS},
\[
(F_\u T_r+F_{\u'}T_r^\#)\cdot \frac12(f_\t\pm f_{\t'})= 
 \frac{\alpha_r(\t)}{\gamma_\u}\frac12(f_\u\pm f_{\u'})
\]
using the fact that $\boldsymbol\alpha$ is an alternating seminormal coefficient system. 
\end{proof}

\begin{Example}
Continuing with Example \ref{spechex1}, noting that $\C$ contains an alternating seminormal coefficient system, and that the generator for $\H_3(\C)^\#$ is $T_1T_2$, using the matrices from Example \ref{spechex1} we obtain two one-dimensional modules on which $T_1T_2$ act by $\omega$ and $\omega^2$ respectively, where $\omega=\frac12(1+\sqrt{-3})\in\C$ is a cube root of unity. 
\end{Example}

Notice that we have the following immediate corollaries of the definitions of the respective modules.   

\begin{Corollary}\label{dimct}
Let $\bla\in\Parts^\ell$ and let $\K$ be a field of characteristic greater than 2. 
\begin{itemize}
\item[(i)] If $\bla\neq\bla'$, $\dim_\K S^{[\bla]}=\dim_\K S^\bla=\left|\Std(\bla)\right|$.  
\item[(ii)] If $\bla=\bla'$, $\dim_\K S^\bla_+=\dim_\K S^\bla_-=\frac12\dim_\K S^\bla=\frac12\left|\Std(\bla)\right|$.  
\end{itemize}
\end{Corollary}

We can now give a classification of irreducible representations for semisimple alternating cyclotomic Hecke algebras. 

\begin{Proposition}\label{spechtheckeballs}
Let $\K$ be a field of characteristic greater than 2 containing an alternating seminormal coefficient system and such that $P_\H$ is nonzero. Then the collection
\[
\{S^{[\bla]}\mid [\bla]\in\(\Parts^\ell)_\sim\mbox{ with } \left|[\bla]\right|=2\}\cup \{S^\bla_+,S^\bla_-\mid [\bla]\in(\Parts^\ell)_\sim\mbox{ with }\left|[\bla]\right|=1\}
\]
is a complete list of pairwise non-isomorphic irreducible modules for the semisimple alternating cyclotomic Hecke algebra $\H_n(\K)^\#$. 
\end{Proposition}
\begin{proof}
By Theorem \ref{ssspecht}, $\sum_{\bla\in\Parts^\ell}\left|\Std(\bla)\right|^2=\ell^nn!$. Moreover, the modules $S^{[\bla]}$ for $\bla\neq\bla'$ and $S^\bla_+$ and $S^\bla_-$ for $\bla=\bla'$ are clearly pairwise non-isomorphic by Lemma \ref{L:separation} and the orthogonality of the idempotents $\{F_\t\mid\t\in\Std(\Parts^\ell)\}$. Hence since
\begin{align*}
\sum_{\substack{[\bla]\in(\Parts^\ell)_\sim\\ \left|[\bla]\right|=2}}\left|\Std(\bla)\right|^2& + \sum_{\substack{[\bla]\in(\Parts^\ell)_\sim\\ \left|[\bla]\right|=1}}2\Bigl(\frac12\left|\Std(\bla)\right|\Bigr)^2\\
&= \sum_{\substack{[\bla]\in(\Parts^\ell)_\sim\\ \left|[\bla]\right|=2}}\left|\Std(\bla)\right|^2+\frac12\sum_{\substack{[\bla]\in(\Parts^\ell)_\sim\\ \left|[\bla]\right|=1}}\left|\Std(\bla)\right|^2\\
&= \frac{\ell^nn!}{2},
\end{align*}
we must have a full list of irreducible representations by Corollary \ref{dim1}, Proposition \ref{dim2} and Proposition \ref{ungradedaltdim}. 
\end{proof}

As a result we obtain the following corollary which again highlights the importance of our seminormal coefficient system. The reader should compare this with the results in \cite{MathasSurvey} and \cite{Riese}. 

\begin{Corollary}
Let $\K$ be a field of characteristic different from 2 containing an alternating seminormal coefficient system and such that $P_\H$ is nonzero. Then $\K$ is a splitting field for $\H_n(\K)^\#$. 
\end{Corollary}
\begin{proof}
The arguments given in this section work over any extension field of $\K$ so the Specht modules constructed in this section are irreducible over any such extension. Hence $\K$ is a splitting field for $\H_n^\#$. 
\end{proof}


    \renewcommand\href[3][\relax]{#3}

\bibliographystyle{andrew}


\end{document}